# Bayesian Learning without Recall

M. Amin Rahimian & Ali Jadbabaie *

*Abstract*—We analyze a model of learning and belief formation in networks in which agents follow Bayes rule yet they do not recall their history of past observations and cannot reason about how other agents' beliefs are formed. They do so by making rational inferences about their observations which include a sequence of independent and identically distributed private signals as well as the actions of their neighboring agents at each time. Successive applications of Bayes rule to the entire history of past observations lead to forebodingly complex inferences: due to lack of knowledge about the global network structure, and unavailability of private observations, as well as third party interactions preceding every decision. Such difficulties make Bayesian updating of beliefs an implausible mechanism for social learning. To address these complexities, we consider a Bayesian without Recall model of inference. On the one hand, this model provides a tractable framework for analyzing the behavior of rational agents in social networks. On the other hand, this model also provides a behavioral foundation for the variety of non-Bayesian update rules in the literature. We present the implications of various choices for the structure of the action space and utility functions for such agents and investigate the properties of learning, convergence, and consensus in special cases.

## I. INTRODUCTION & BACKGROUND

Individuals often exchange opinions with their peers in order to learn from their knowledge and experiences, and in making various decisions such as investing in stock markets, voting in elections, choosing their political affiliations, selecting a brand of a product or a medical treatment. These interactions occur through a variety of media which we collectively refer to as social networks. James Surowiecki in his popular science book on wisdom of crowds [1], provides well-known cases for information aggregation in social networks, and argues how under the right circumstances (diversity of opinion, independence, decentralization and aggregation) groups outperform even their smartest or best informed members; see for example the essentially perfect performance of the middlemost estimate at the weight-judging competition of the 1906 West of England Fat Stock and Poultry Exhibition studied by Francis Galton in his 1907 Nature article [2], entitled "Vox Populi" (The Wisdom of Crowds), or the study of market reaction to the 1986 challenger disaster in [3], where its is pointed out that the main responsible company's (Morton Thiokol) stock was hit hardest of all, even months before the cause of the accident could be officially determined.

On the other hand, several studies point out that the evolution of people's opinions and decisions in social networks is subject to various kind of biases and inefficiencies [4]–[8]. Such deviations from the efficient and/or rational outcome are often attributed to the structural effects that arise in networked interactions; in particular, the predominant influence of more central agents in shaping the group decision, in spite of the fact that such influential agents do not necessarily enjoy a high quality of observations or superior knowledge; cf. persuasion bias in [9], obstructions to wisdom of crowds in [10], and data incest in [11]. Subsequently, a better understanding of social learning can also help us analyze the effect of such biases and use our insights and conclusions to improve policy designs that are aimed at implementing desirable social norms or eradicating undesirable ones, or even to come up with more efficient procedures for aggregating individual beliefs, and to understand how media sources, prominent agents, government and politicians are able to manipulate public opinion and influence spread of beliefs in society [12].

We model the set of possible alternatives that are of common interest to all individual in society by a set of finitely many states of the world; and endow each individual agent with a belief, representing her opinion or understanding of the true state of the world. Thereby, agents exchange beliefs in social networks to benefit from each other's opinions and private information in trying to learn an unknown state of the world. The problem of social learning is to characterize and understand such interactions and it is a classical focus of research in behavioral microeconomic theory [13, Chapter 8], [14, Chapter 5], [15]. Research on formation and evolution of beliefs in social networks and subsequent shaping of the individual and mass behavior has attracted much attention amongst diverse communities in engineering [16]–[18], statistics [19], economics [20], and sociology [21]. The problem of social learning has close siblings in distributed estimation [22], [23], data fusion [24], and statistical learning theory [25]; while relations can be also traced to the consensus and coordination problems that are studied in the distributed control theory [26], [27].

Consider an agent trying to estimate an unknown state of the world. She bases her estimation on a sequence of independent and identically distributed (i.i.d.) private signals that she observes and whose common distribution is determined by the unknown state. Suppose further that her belief about the unknown state is represented by a discrete probability distribution over the set of finitely many possibilities $\Theta$, and that she sequentially applies Bayes rule to her observations at each step, and updates her beliefs accordingly. It is a well-known consequence of the classical results in merging and learning theory [28], [29] that the beliefs formed in the above manner constitute a bounded martingale and converge to a limiting distribution as the number of observations increases. However, the limiting distribution may differ from a point mass centered at the truth, in which case the agent fails to learn the true state asymptotically. This may be the case, for instance if the agent faces an identification problem, that is when there are states other than the true state which are observationally equivalent to the true state and induce the same distribution on her sequence of privately observed signals. Accordingly, the agents have an incentive to communicate in a social network so that they can resolve their identification problems by relying on each other's observational abilities.

*Correspondence to: Ali Jadbabaie, Institute for Data, Systems, and Society (IDSS), Massachusetts Institute of Technology (MIT), Cambridge, MA 02139, USA. (email: jadbabai@mit.edu). This work was supported by ARO MURI W911NF-12-1-0509.

Rational agents in a social network would apply Bayes rule successively to their observations at each step, which include not only their private signals but also the beliefs and actions communicated by their neighbors. However, such repeated applications of Bayes rule in networks become very complex, especially if the agents are unaware of the global network structure. This is due to the fact that the agents at each step should use their local data that is increasing with time, and make very complex inferences about possible signal structures leading to their observations. Indeed, tractable modeling and analysis of rational behavior in networks is an important problem in network economics and have attracted much attention, [19], [30], [31].

Modeling memory constraints in the context of social learning is an important research problem [15, Chapter 5]. In recent results, Wilson [32] considers the model of a decision maker who chooses between two actions with pay-offs that depend on the true state of the world. Furthermore, the decision maker must always summarize her information into one of the finitely many states, leading to optimal decision rules that specify the transfers between these states. The problem of learning with finite memory in the context of hypothesis testing was originally formulated by [33], [34], where memory constraints restrict the storage capacity for the test statistics. Accordingly, while sufficient statistics are very useful computational tools their utility for memory reduction is not clear. Subsequent results provide sophisticated algorithms to perform the task of hypothesis testing using test statistics that take only finitely many values and to guarantee an asymptotically vanishing error probability [35]–[38]. More recently, the authors in [39] have considered this problem in a setting where agents each receive an independent private signal and make decisions sequentially. Memory in this context refers to the number of immediate predecessors whose decisions are observable by any given agent at the time of making her decision. Accordingly, while the almost sure convergence of the sequence of individual decisions to the correct state is not possible in this finite memory setting, the authors construct decision rules that achieve convergence and learning in probability. They next go on to consider the behavior of rational (pay-off maximizing) agents in this context and show that in no equilibrium of the associated Bayesian game learning can occur.

To avoid the complexities of fully rational inference, a variety of non-Bayesian update rules have been proposed that rely on the seminal work of DeGroot in linear opinion pooling [40], where agents update their opinions to a convex combination of their neighbors' beliefs and the coefficients correspond to the level of confidence that each agent puts in each of her neighbors. More recently, [20], [41] consider a variation of this model for streaming observations, where in addition to the neighboring beliefs the agents also receive private signals. Other forms of non-Bayesian rules were studied by [42] who consider a variation of observational learning in which agents observe the action and pay-offs of their neighbors and make rational inferences about these action/pay-off correspondence together with the choices made by their neighbors, but ignore the fact that their neighbors are themselves learning from their own observations. In more recent results, [43], [44] consider models of autarkic play where players at each generation observe their predecessor but naïvely think that any predecessor's action relies solely on that player's private information, thus ignoring the possibility that successive generations are learning from each other.

*A. Motivation & Contributions*

In this paper, we analyze a model of repeated interactions for social learning when agents receive private signals and observe their neighboring decisions (actions) at every epoch of time. Such a model is a good descriptor for online reputation and polling systems such as Yelp® and TripAdvisor®, where individuals' recommendations are based on their private observations and recommendations of their friends [45, Chapter 5]. The analysis of such systems is important not only because they play a significant role in generating revenues for the businesses that are being ranked [46], but also for the purposes of designing fair rankings and accurate recommendation systems.

Heuristics are widely used in the literature to model social interactions and decision making [47]–[49]. They provide tractable tools to analyze boundedly rational behavior and offer useful insights about decision making under uncertainty [10], [50]. They are also verified to be good descriptors for the behavior of real world agents in the experimental studies by Grimm and Mengel [51] and Chandrasekhar, Larreguy and Xandri [52]. Despite their widespread applications, theoretical and axiomatic foundations of social inferences using non-Bayesian (heuristic) update rules have received limited attention and only recently [53], [54]. A comprehensive theory of non-Bayesian learning that reconciles the rational and boundedly rational approaches with the widely used heuristics remains in demand. A chief contribution of this paper is in establishing a behavioral foundation for the existing non-Bayesian updates in the literature. In particular, this paper addresses the question of how one can limit the information and cognitive requirements of a Bayesian inference and still ensure consensus or learning for the agents.

Some of the non-Bayesian update rules have the property that they resemble the replication of a first step of a Bayesian update from a common prior cf. [9]; in this paper, we formalize and expand up on this idea. In particular, we propose the so-called *Bayesian without Recall* (BWR) model as a belief formation and update rule for *Rational but Memoryless* agents.

To model the behavior of such agents we replicate the time-one update of a Bayesian agent for all future time steps. On the one hand, the BWR model of inference is motivated by the real-world behavior of people reflected in their spur-of-the-moment decisions and impromptu behavior; basing decisions only on the immediately observed actions and without regard for the history of such actions. On the other hand, BWR offers a boundedly rational approach to model decision making over social networks. The latter is in contrast with the Bayesian approach which is not only unrealistic in the amount of cognitive burden that it imposes on the agents, but also is often computationally intractable and complex to analyze.

*B. Brief Overview & Organization of Paper*

A key message of this paper is to show how the BWR scheme can provide justification for some well-known update rules. These rules can be equivalently explained as the update of a Bayesian agent that naively assumes the beliefs of each of her neighbors were formed by a private observation, and not through repeated interaction with others. We begin by specializing the BWR model to a case where agents try to decide between one of the two possible states and are rewarded for every correct choice that they make. We show that the BWR action updates in this case are given by *weighted majority and threshold rules* that linearly combine the observed binary actions and log-likelihoods of private signal. We show that under these update rules the action profiles evolve as a Markov chain on the Boolean cube and that the properties of consensus and learning are subsequently determined by the equilibria of this Markov chain.

When there are only finitely many states of the world and agents choose actions over the probability simplex, then the action spaces are rich enough to reveal the beliefs of every communicating agent. We show that the BWR updates in this second case are *log-linear* in the reported beliefs of the neighbors and the likelihood of private signals. We investigate the properties of convergence and learning for such agents in a strongly connected social network, provided that the truth is identifiable through the aggregate observations of the agents across entire network. This is of particular interest, when the agents cannot distinguish the truth based solely on their private observations, and yet together they learn. Analysis of convergence and learning in this case reveals that almost-sure learning happens only if the agents are arranged in a directed circle. We explain how the circular BWR updates generalize to any strongly connected topology by choosing a single neighbor randomly at every time step. We characterize the rate of learning in such cases as being asymptotically exponentially fast with an exponent that is linear in time and whose coefficient can be expressed as a weighted average of the relative entropies of the signal likelihoods of all agents.

The remainder of this paper is organized as follows. The details of our BWR model for social learning are explained in Section II. In Section III we specialize the BWR model to a binary state and action space and investigate the evolution of actions in the resultant Ising model. Next in Section IV we analyze the case where the network agents are announcing their beliefs at every epoch and the BWR updates become log-linear. Section V concludes the paper, followed by four mathematical appendices, labeled A to E, which provide the detailed derivations and proofs for the technical results that are presented throughout the paper.

## II. THE BWR MODEL OF NAIVE INFERENCE

A dual process theory for the psychology of thinking identifies two systems for the operations of mind [55]: one that is fast, intuitive, non-deliberative, habitual and automatic (system one); and a second one that is slow, attentive, effortful, deliberative, and conscious (system two) [56]. Major advances in behavioral economics are due to incorporation of this dual process theory and the subsequent models of bounded rationality [57]. Reliance on heuristics for decision making is a distinctive feature of system one that avoids the computational burdens of a rational evaluation; system two on the other hand, is bound to deliberate on the options based on the available information before making recommendations. The interplay between these two systems and how they shape the individual decisions is of paramount importance [58]. According to the BWR model of naive inference for social learning, as the agent experiences with her environment her initial response would engage her system two: she rationally evaluates the reports of her neighbors and use them along with her own private signal to make a decision. However, after her initial experience and by engaging in repeated interactions with the environment, her system one takes over her decision processes, implementing a heuristic that imitates her (rational/Bayesian) inferences from her initial experience; hence avoiding the burden of additional cognitive processing in the ensuing interactions with her environment. In the following subsections we set forth the mathematical and modeling details that allow us to implement the BWR model in a variety of environments with different utility, observation and information structures.

*A. Signals and Environment*

Consider a set of $n$ agents that are labeled by $[n]$ and interact according to a digraph $\mathcal{G} = ([n], \mathcal{E})$.[1] The neighborhood of agent $i$ is the set of all agents who communicate with agent $i$ and it is denoted by $\mathcal{N}(i) = \{j \in [n]; (j,i) \in \mathcal{E}\}$; the degree of agent $i$ is the cardinality of its neighborhood $\text{card}(\mathcal{N}(i))$. Digraph $\mathcal{G}$ is strongly connected if there are directed paths from any agent to any other agents.

Let $\Theta$ denote a finite set of possible states of the world and $\Delta\Theta$ be the space of all probability measures on the set $\Theta$.

---

[1] Throughout the paper, $\mathbb{R}$ is the set of real numbers, $\mathbb{N}$ denotes the set of all natural numbers, and $\mathbb{N}_0 := \mathbb{N} \cup \{0\}$. For $n \in \mathbb{N}$ a fixed integer the set of integers $\{1, 2, \ldots, n\}$ is denoted by $[n]$, while any other set is represented by a calligraphic capital letter. The cardinality of a set $\mathcal{X}$, which is the number of its elements, is denoted by $\text{card}(\mathcal{X})$. The set difference between $\mathcal{X}$ and $\mathcal{Y}$ denoted by $\mathcal{X}\setminus\mathcal{Y}$ is the set of all elements of $\mathcal{X}$ that do not belong to $\mathcal{Y}$.

A parameter $\theta \hookleftarrow \Theta$ is chosen arbitrarily from $\Theta$ by nature.[1] Associated with each agent $i$, $\mathcal{S}_i$ is a finite set called the signal space of $i$ and given $\theta$, $l_i(\cdot|\theta)$ is a probability mass function on $\mathcal{S}_i$, which is referred to as the *signal structure* or *likelihood function* of agent $i$. Let $t \in \mathbb{N}_0$ denote the time index; for each agent $i$, $\{\mathbf{s}_{i,t}, t \in \mathbb{N}_0\}$ is a sequence of independent and identically distributed (i.i.d.) random variables that take values in $\mathcal{S}_i$ and with the probability mass function $l_i(\cdot|\theta)$. This sequence represents the private signals that agent $i$ observes over time. Note that the private signals are independent over time and across the agents.

**Remark 1** (Learning from others). *The fact that different people make independent observations about the underlying truth state $\theta$ gives them incentive to communicate in social networks, in order to benefit from each others' observations and to augment their private information. Moreover, different people differ in their observational abilities. For instance suppose that the signal structure of agent $i$ allows her to distinguish the truth $\theta$ and the false state $\check{\theta}$, while the two states $\hat{\theta}$ and $\theta$ are indistinguishable to her: i.e. $\ell_i(s_i|\check{\theta}) \neq \ell_i(s_i|\theta)$ for some $s_i \in \mathcal{S}_i$, whereas $\ell_i(s_i|\hat{\theta}) = \ell_i(s_i|\theta)$ for all $s_i \in \mathcal{S}_i$. In such circumstances, agent $i$ can never resolve her ambiguity between $\theta$ and $\hat{\theta}$ on her own; hence, she has no choice but to rely on other people's observations to be able to learn the truth state with certainty.*

For any $\check{\theta} \in \Theta$ let $\lambda_{\check{\theta}} : \cup_{i \in [n]} \mathcal{S}_i \to \mathbb{R}$ be the real valued function measuring log-likelihood ratio of the signal $s_i$ under states $\check{\theta}$ and $\theta$, defined as $\lambda_{\check{\theta}}(s_i) := \log \left( \ell_i(s_i|\check{\theta})/\ell_i(s_i|\theta) \right)$. This is a measure of the information content that the signal $s_i$ provides for distinguishing the false state $\check{\theta}$ from the truth $\theta$. Subsequently we work with the probability triplet $(\Omega, \mathscr{F}, \mathbb{P}_\theta)$, where $\Omega = \left( \prod_{i \in [n]} \mathcal{S}_i \right)^{\mathbb{N}_0}$ is an infinite product space with a typical element $\omega = ((s_{1,0}, \ldots, s_{n,0}), (s_{1,1}, \ldots, s_{n,1}), \ldots)$ and the associated sigma field $\mathscr{F} = \mathscr{P}(\Omega)$. The probability measure on $\Omega$ is $\mathbb{P}_\theta(\cdot)$ which assigns probabilities consistently with the likelihood functions $l_i(\cdot|\theta), i \in [n]$; and in such a way that conditional on $\theta$ the random variables $\mathbf{s}_{i,t}, i \in [n], t \in \mathbb{N}_0$ taking values in $\mathcal{S}_i, i \in [n]$, are independent. The expectation operator $\mathbb{E}_\theta\{\cdot\}$ represents integration with respect to $d\mathbb{P}_\theta(\omega), \omega \in \Omega$.

### B. Beliefs, Observations, Actions and Rewards

An agents' belief about the unknown allows her to make decisions even as her pay-off is dependent on the unknown state $\theta$. These beliefs about the unknown state are probability distributions over $\Theta$. Even before the true state $\theta$ is assigned and any observations are made, every agent $i \in [n]$ holds a prior belief $\nu_i(\cdot) \in \Delta\Theta$ with full support: $\nu_i(\hat{\theta}) > 0$, $\forall \hat{\theta} \in \Theta$; this represents her subjective biases about the true value of $\theta$. For each time instant $t$, let $\boldsymbol{\mu}_{i,t}(\cdot)$ be probability mass function on $\Theta$, representing the *opinion* or *belief* at time $t$ of agent $i$ about the realized value of $\theta$, and define $\boldsymbol{\phi}_{i,t}(\check{\theta}) := \log \left( \boldsymbol{\mu}_{i,t}(\check{\theta})/\boldsymbol{\mu}_{i,t}(\theta) \right)$ as the log-belief ratio of agent $i$ at time $t$ under the states $\check{\theta}$ and $\theta$. Moreover, let $\mathbb{P}_{i,t}\{\cdot\}$ denote the probability measure on $\Theta \times \Omega$ that assigns probabilities consistently with $\boldsymbol{\mu}_{i,t}(\cdot)$ and the independent signals likelihoods, and let the associated expectation operator be $\mathbb{E}_{i,t}\{\cdot\}$.

At $t = 0$ after $\theta \hookleftarrow \Theta$ is assigned, the values $s_i \in \mathcal{S}_i$ of $\mathbf{s}_{i,0}$ are realized and the latter is observed privately by each agent $i$ for all $i \in [n]$. Associated with every agent $i$ is an action space $\mathcal{A}_i$ that represents all the choices available to her at every point of time $t \in \mathbb{N}_0$, and a utility function $u_i(\cdot, \cdot) : \mathcal{A}_i \times \Theta \to \mathbb{R}$. Subsequently, at every time $t \in \mathbb{N}_0$ each agent $i \in [n]$ chooses an action $\mathbf{a}_{i,t} \in \mathcal{A}_i$ and is rewarded $u_i(\mathbf{a}_{i,t}, \theta)$. [2]

Such modeling of rewards to actions at successive time periods is common place in the study of learning in games and a central question of interest is that of regret which measures possible gains by the players if they were to play other actions from what they have chosen in a realized path of play; in particular, existence of update/decision rules that would guarantee a vanishing time-average of regret as $t \to \infty$, and possible equilibria that characterize the limiting behavior of agents under such rules have attracted much attention, cf. [25], [60], [61]. Under a similar framework Rosenberg et al. [62] consider consensus in a general setting with players who observe a private signal, choose an action and receive a pay-off at every stage, and pay-offs that depend only on an unknown parameter and players' actions. They show that in this setting with no pay-off externalities and interactions which are purely informational players asymptoticly play their best-replies given their beliefs and will agree in their pay-offs; in particular, all motives for experimentation will eventually disappear.

### C. Naive (Memoryless) Action Updates

Given $\mathbf{s}_{i,0}$, agent $i$ forms an initial Bayesian opinion $\boldsymbol{\mu}_{i,0}(\cdot)$ about the value of $\theta$, which is given by

$$\boldsymbol{\mu}_{i,0}(\hat{\theta}) = \frac{\nu_i(\hat{\theta}) l_i(\mathbf{s}_{i,0} \mid \hat{\theta})}{\sum_{\tilde{\theta} \in \Theta} \nu_i(\tilde{\theta}) l_i(\mathbf{s}_{i,0} \mid \tilde{\theta})}, \forall \hat{\theta} \in \Theta. \quad (1)$$

She then chooses the action: $\mathbf{a}_{i,0} \hookleftarrow \arg\max_{a_i \in \mathcal{A}_i} \sum_{\hat{\theta} \in \Theta} u_i(a_i, \hat{\theta}) \boldsymbol{\mu}_{i,0}(\hat{\theta})$, maximizing her expected reward: $\mathbb{E}_{i,0}\{u_i(\mathbf{a}_{i,0}, \theta)\}$. Not being notified of the actual realized value for $u_i(\mathbf{a}_{i,0}, \theta)$, she then observes the actions that her neighbors have taken: $\mathbf{a}_{j,0}, j \in \mathcal{N}(i)$. Given her extended set of observations at time $t = 1$, she makes a second and possibly different move $\mathbf{a}_{i,1}$ according to

$$\mathbf{a}_{i,1} \hookleftarrow \arg\max_{a_i \in \mathcal{A}_i} \sum_{\hat{\theta} \in \Theta} u_i(a_i, \hat{\theta}) \boldsymbol{\mu}_{i,1}(\hat{\theta}), \quad (2)$$

maximizing her expected pay off conditional on everything that she has observed thus far: $\mathbb{E}_{i,1}\{u_i(\mathbf{a}_{i,1}, \hat{\theta})\} = \mathbb{E}\{u_i(\mathbf{a}_{i,1}, \hat{\theta}) \mid \mathbf{s}_{i,0}, \mathbf{a}_{j,0} : j \in \mathcal{N}(i)\}$. Subsequently, she is granted her net reward of $u_i(\mathbf{a}_{i,0}, \theta) + u_i(\mathbf{a}_{i,1}, \theta)$ from her past two plays.

---

[1] For a set $\mathcal{X}$, $x \hookleftarrow \mathcal{X}$ denotes an arbitrary choice from the elements of $\mathcal{X}$ that is assigned to $x$. The power-set of $\mathcal{X}$ is the set of all its subsets and it is denoted by $\mathscr{P}(\mathcal{X}) = \{\mathcal{M}; \mathcal{M} \subset \mathcal{X}\}$. Boldface letters denote random variables, vectors are denoted by a bar over their respective upper or lower-case letters and $^T$ denotes matrix transpose.

[2] The utility functions $u_i(\cdot, \cdot)$, signal structures $l_i(\cdot|\cdot)$, priors $\nu_i(\cdot)$, as well as the corresponding sample spaces $\mathcal{A}_i, \mathcal{S}_i$ and $\Theta$ are all common knowledge amongst the communicating agents for all $i \in [n]$. The assumption of common knowledge in the case of fully rational (Bayesian) agents implies that given the same observations of one another's actions or private signals distinct agents would make identical inferences; in the sense that starting form the same belief about the unknown $\theta$, their updated beliefs given the same observations would be the same; in Aumann's words, rational agents cannot agree to disagree [59].

Following realization of rewards for their first two plays, in any subsequent time instance $t > 1$ each agent $i \in [n]$ observes a private signal $\mathbf{s}_{i,t}$ together with the preceding actions of her neighbors $\mathbf{a}_{j,t-1}$, $j \in \mathcal{N}(i)$. She then takes an option $\mathbf{a}_{i,t}$ out of the set $\mathcal{A}_i$, such that her expected utility given her observations is maximized. Of particular significance in our description of the behavior of agents in the succeeding time periods $t > 1$, is the relation

$$f_i(\mathbf{s}_{i,0}, \mathbf{a}_{j,0} : j \in \mathcal{N}(i)) := \mathbf{a}_{i,1} \hookleftarrow \arg\max_{a_i \in \mathcal{A}_i} \mathbb{E}_{i,1}\{u_i(a_i, \hat{\theta})\} \quad (3)$$

derived in (2), which given the observations of agent $i$ from time $t = 0$, specifies her (Bayesian) pay-off maximizing action for time $t = 1$. Note that in writing (2), we assumed that the agents do not receive any private signals at $t = 1$ and there is therefore no $\mathbf{s}_{i,1}$ appearing in the updates of any agent $i$; and this convention is exactly to facilitate the derivation of mapping $f_i : \mathcal{S}_i \times \prod_{j \in \mathcal{N}(i)} \mathcal{A}_j \to \mathcal{A}_i$, from the private signal space and action spaces of the neighbors to succeeding actions of each agent. In every following instance we aim to model the inferences of agents about their observations as being rational but memoryless: as of those who come to know their immediate observations which include the actions of their neighbors and their last private signals, but cannot trace these observations to their roots and has no ability to reason about why their neighbors may be behaving the way they do. In particular, such agents have no incentives for experimenting with false reports, as their lack of memory prevents them from reaping the benefits of their experiment, including any possible revelations that a truthful report may not reveal. Subsequently, we argue on normative grounds that such rational but memoryless agents would replicate the behavior of a Bayesian (fully-rational) agent between times zero and one; whence by regarding their observations as being direct consequences of inferences that are made based on the initial priors, they reject any possibility of a past history beyond their immediate observations:[1]

$$\mathbf{a}_{i,t} = f_i(\mathbf{s}_{i,t}, \mathbf{a}_{j,t-1} : j \in \mathcal{N}(i)), \forall t > 1.$$

On the other hand, note that rationality of agents constrains their beliefs $\boldsymbol{\mu}_{i,t}(\cdot)$ given their immediate observations; hence, we can also write

$$\mathbf{a}_{i,t} \hookleftarrow \arg\max_{a_i \in \mathcal{A}_i} \mathbb{E}_{i,t}\{u_i(a_i, \hat{\theta})\}, \quad (4)$$

or equivalently $\mathbf{a}_{i,t} \hookleftarrow \arg\max_{a_i \in \mathcal{A}_i} \sum_{\hat{\theta} \in \Theta} u_i(a_i, \hat{\theta}) \boldsymbol{\mu}_{i,t}(\hat{\theta})$.

In the sequel, we explore various structures for the action space and the resultant update rules $f_i$. In Section III, we show how a common heuristic such as weighted majority can be explained as a rational but memoryless behavior with actions taken from a binary set. In Section IV we shift focus to a finite state space and the probability simplex as the action space. There agents exchange beliefs and the belief updates are log-linear.

### III. WEIGHTED MAJORITY AND THRESHOLD RULES

Consider a binary state space $\Theta = \{+1, -1\}$, and suppose that the agents have a common binary action space $\mathcal{A}_i = \{-1, 1\}$, for all $i$. Let their utilities be given by $u_i(a, \theta) = 2\mathbb{1}_a(\theta) - 1$, for any agent $i$ and all $\theta, a \in \{-1, 1\}$; here, $\mathbb{1}_a(\theta)$ is equal to one only if $\theta = a$ and is equal to zero otherwise. Subsequently, the agent is rewarded by $+1$ every time she correctly determines the value of $\theta$ and is penalized by $-1$ otherwise (Fig. 1).

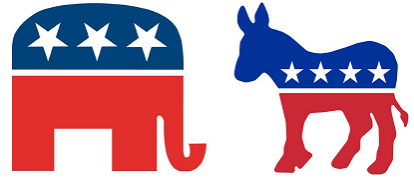

Fig. 1: Bipartisanship is an example of a binary state space.

We can now calculate $\sum_{\hat{\theta} \in \Theta} u_i(a_i, \hat{\theta}) \boldsymbol{\mu}_{i,t}(\hat{\theta}) = a(\boldsymbol{\mu}_{i,t}(+1) - \boldsymbol{\mu}_{i,t}(-1)) = a(2\boldsymbol{\mu}_{i,t}(+1) - 1), \forall a \in \{-1, 1\}$;

---

[1]Extensions of the above behavioral model to rational agents with bounded memory is of interest; nonetheless, the analysis of Bayesian update even in the simplest cases become increasingly complex. As an example, consider a rational agent who recalls only the last two epochs of her past. In order for such an agent to interpret her observations in the penultimate and ultimate steps, she needs not only a prior to interpret her neighbor's beliefs at the penultimate step, but also a prior to interpret her neighbor's inferences about what she reported to them at the penultimate step leading to their ultimate beliefs. In other words, she needs a prior on what her neighbor's regard as her prior when they interpret what she reports to them as her penultimate belief. Indeed, such belief hierarchies are commonplace in game-theoretic analysis of incomplete information and are captured by the formalism of type space [63], [64].

and from (4), we get[1]

$$\mathbf{a}_{i,t} = \begin{cases} 1 & \text{if } \boldsymbol{\mu}_{i,t}(+1) \geq \boldsymbol{\mu}_{i,t}(-1), \\ -1 & \text{if } \boldsymbol{\mu}_{i,t}(+1) < \boldsymbol{\mu}_{i,t}(-1), \end{cases} \quad (5)$$

We can now proceed to derive the memoryless update rule $f_i$ under the above prescribed settings. This is achieved by the following expression of the action update of agent $i$ at time 1. Throughout this section and without any loss of generality, we assume that $\theta = -1$.

**Lemma 1** (Time-One Bayesian Actions). *The Bayesian action of agent $i$ at time one following her observations of actions of her neighbors at time zero and her own private signal at time zero is given by $\mathbf{a}_{i,1} = \text{sign}(\sum_{j \in \mathcal{N}(i)} w_j \mathbf{a}_{j,0} + \eta_i + \lambda_1(\mathbf{s}_{i,0}))$, where $w_i$ and $\eta_i$ are constants for each $i$ and they are completely determined by the initial prior and signal structures of agent $i$ and her neighbors.*

The exact expressions of the constants $w_i$, $\eta_i$ and their derivations can be found in Appendix A. Indeed, making the necessary substitutions we derive the following memoryless update $f_i$ for all $t > 1$: $\mathbf{a}_{i,t} = \text{sign}\left(\sum_{j \in \mathcal{N}(i)} w_j \mathbf{a}_{j,t-1} + \eta_i + \lambda_1(\mathbf{s}_{i,t})\right)$. This update rule has a familiar format as a weighted majority and threshold function with the weights and threshold given by $w_i$ and $\mathbf{t}_{i,t} := -\lambda_1(\mathbf{s}_{i,t}) - \eta_i$, the latter being random and time-varying.

Majority and threshold functions are studied in the analysis of Boolean functions [65, Chapter 5] and several properties of them including their noise stability are of particular interest [66]–[68]. This update rule also appears as the McCulloch-Pitts model of an artificial neuron [69], with important applications in neural networks and computing [70]. This update rule is also important in the study of the Glauber Dynamics in the Ising model, where the $\pm 1$ states represent atomic spins. The spins are arranged in a graph and each spin configuration has a probability associated with it depending on the temperature and the interaction structure [71, Chapter 15], [72]. The Ising model provides a natural setting for the study of cooperative behavior in social networks. Recent studies have explored the applications of Ising model for analysis of social and economic phenomena such as rumor spreading [73], study of market equilibria [74], and opinion dynamics [75].

Following this model, every agent $i \in [n]$ chooses her action $\mathbf{a}_{i,t} \in \{\pm 1\}$ as the sign of $\sum_{j \in \mathcal{N}(i)} w_j \mathbf{a}_{j,t-1} + \eta_i + \lambda_1(\mathbf{s}_{i,t})$. Subsequently, in processing her available data and choosing her action $\mathbf{a}_{i,t}$, every agent seeks to maximize $(\sum_{j \in \mathcal{N}(i)} w_j \mathbf{a}_{j,t-1} \mathbf{a}_{i,t} + \mathbf{a}_{i,t}(\eta_i + \lambda_1(\mathbf{s}_{i,t})))$. Hence, we can interpret each of the terms appearing as the argument of the sign function, in accordance with how they influence agent $i$'s choice of action. In particular, the term $\eta_i + \lambda_1(\mathbf{s}_{i,t})$ represents the propensity of agent $i$ in choosing the false action $\theta_1 := 1$ at time $t$, and it is determined by the log-likelihood ratio of

---
[1]In writing (5) we follow the convention that agents choose +1 when they are indifferent between their two options. Similarly, the sign function is assumed to take the value +1 when its argument is zero. This assumption is consistently followed everywhere throughout this paper, except in Proposition 1 and its proof in Appendix C, see the footnote therein for further details.

---

private signal $\lambda_1(\mathbf{s}_{i,t})$, as well as her innate tendency towards +1 irrespective of any observations. The latter is reflected in the constant $\eta_i := \log(\nu_i(\theta_1)/\nu_i(\theta_2)) + \log V_i$ based on the log-ratio of her initial prior belief and her knowledge of her neighbor's signal structures, as captured by the constant $V_i$ in (17) of Appendix A. The latter is increasing in $\ell_j(s_j \mid \theta_1)$ and decreasing in $\ell_j(s_j \mid \theta_2)$ for any fixed signal $s_j \in \mathcal{S}_j$, $j \in \mathcal{N}(i)$; cf. Lemma 3 of Appendix A.

By the same token, we can also interpret the interaction terms $w_j \mathbf{a}_{j,t-1} \mathbf{a}_{i,t}$. Lemma 4 of Appendix A establishes that constants $w_j$ are non-negative for every agent $j \in [n]$. Hence, in maximizing $\sum_{j \in \mathcal{N}(i)} w_j \mathbf{a}_{j,t-1} a + a(\eta_i + \lambda_1(\mathbf{s}_{i,t}))$ through her choice of $a \in \pm 1$ at every time $t$, agent $i$ aspires to align her choice with as many of her neighbors $j \in \mathcal{N}(i)$ as possible. However, in doing so she weighs more the actions of those of her neighbors $j \in \mathcal{N}(i)$ who have larger constants $w_j$. The constant $w_j := \log W_j$ with $W_j$ given in (18) of Appendix A is a measure of observational ability of agent $j$ as relates to our model: agents with large constants $w_j$ are those who hold expert opinions in the social network and they play a major role in shaping the actions of their neighboring agents. Positivity of $w_i$ for any $i \in [n]$, per Lemma 4 of Appendix A, also signifies a case of positive externalities: an agent is more likely to choose an action if her neighbors make the same decision.

### A. Analysis of Convergence and Learning in Ising Networks

To begin with the analysis of the binary action update dynamics derived above, we introduce some useful notation. For all $t \in \mathbb{N}_0$, let $\overline{\mathbf{a}}_t := (\mathbf{a}_{1,t}, \ldots, \mathbf{a}_{n,t})^T$ be the profile of actions taken by all agents at time $t$. Subsequently, we are interested in the probabilistic evolution of the action profiles $\overline{\mathbf{a}}_t, t \in \mathbb{N}_0$ under the following dynamics

$$\mathbf{a}_{i,0} = \text{sign}\left(\log \frac{\nu_i(\theta_1)}{\nu_i(\theta_2)} + \lambda_1(\mathbf{s}_{i,0})\right), \quad (6)$$

$$\mathbf{a}_{i,t} = \text{sign}\left(\sum_{j \in \mathcal{N}(i)} w_j \mathbf{a}_{j,t-1} + \eta_i + \lambda_1(\mathbf{s}_{i,t})\right), t \geq 1, \quad (7)$$

for all $i \in [n]$. The two constants $w_i$ and $\eta_i$ for each agent $i$ are specified in Appendix A and they depend only on the signal structure and initial prior of that agent and her neighbors. The evolution of action profiles $\overline{\mathbf{a}}_t$ in (7) specifies a finite Markov chain that jumps between the vertices of the Boolean hyper cube, $\{\pm 1\}^n$. The analysis of the time-evolution of action profiles is facilitated by the classical results from the theory of finite Markov chains with the details spelled out in Appendix B.

If the signal structures are rich enough to allow for sufficiently strong signals (having large absolute log-likelihood ratios), or if the initial priors are sufficiently balanced (dividing the probability mass almost equally between $\theta_1$ and $\theta_2$), then any action profiles belonging to $\{\pm 1\}^n$ is realizable as $\overline{\mathbf{a}}_0$ with positive probability under (6). In particular, any recurrent state of the finite Markov chain over the Boolean cube is reachable with positive probability and the asymptotic behavior can be only determined up to a distribution over the first set of communicating recurrent states that is reached by $\bar{\mathbf{a}}_t$, cf. Proposition 2 of Appendix B. However, if a recurrent class constitutes a singleton, then our model makes sharper predictions: $\lim_{t\to\infty} \bar{\mathbf{a}}_t$ almost surely exists and is identified as an absorbing state of the finite Markov chain. This special case is treated next due to its interesting implications.

*B. Equilibrium, Consensus, and (Mis-)Learning*

We begin by noting that the absorbing states of the Markov chain of action profiles specify the equilibria under the action update dynamics in (7). Formally, an equilibrium $\bar{a}^* \in \{\pm 1\}^n$ is such that if the dynamics in (7) is initialized by $\bar{\mathbf{a}}_0 = \bar{a}^*$, then with probability one it satisfies $\bar{\mathbf{a}}_t = \bar{a}^*$ for all $t \geq 1$. Subsequently, the set of all equilibria is completely characterized as the set of all absorbing states, i.e. any action profiles $\bar{a}^* \in \{\pm 1\}^n$ satisfying $P(\bar{a}^*, \bar{a}^*) = 1$, where $P : \{\pm 1\}^n \times \{\pm 1\}^n \to [0,1]$ specifies the transition probabilities in the Markov chain of action profiles, as defined in (21) of Appendix B. It is useful to express this condition in terms of the model parameters as follows. The proof is included in Appendix C and with a caveat explained in its footnote.

**Proposition 1** (Characterization of the Equilibria). *An action profile $(a_1^*, \ldots, a_n^*) \in \{\pm 1\}^n$ is an equilibrium of (7) if, and only if, $-\min_{s_i \in \mathcal{S}_i} a_i^*(\lambda_1(s_i) + \eta_i) \leq \sum_{j \in \mathcal{N}(i)} w_j a_j^* a_i^*$, $\forall i \in [n]$.*

Of particular interest are the two action profiles $(1, \ldots, 1)^T$ and $(-1, \ldots, -1)^T$ which specify a consensus amongst the agents in their chosen actions. The preceding characterization of equilibria is specialized next to highlight the necessary and sufficient conditions for the agents to be at equilibrium whenever they are in consensus.

**Corollary 1** (Equilibrium at Consensus). *The agents will be in equilibrium at consensus if, and only if, $\max_{s_i \in \mathcal{S}_i} |\lambda_1(s_i) + \eta_i| < \sum_{j \in \mathcal{N}(i)} w_j$, $\forall i \in [n]$.*

The requirement of learning under our model is for the agents to reach a consensus on truth. That is for the action profiles $\bar{\mathbf{a}}_t$ to converge to $(\theta, \ldots, \theta)$ as $t \to \infty$. In particular, as in Corollary 1, we need agents to be at equilibrium when in consensus; hence, there would always be a positive probability for the agents to reach consensus on an untruth: with a positive probability, the agents (mis-)learn.

Next in Section IV we show that when the action space is rich enough to reveal the beliefs of the agents, then the rational but memoryless behavior culminates in a log-linear updating of the beliefs with the observations. The analysis of convergence and learning under these log-linear updates consumes the bulk of that section.

## IV. LOG-LINEAR UPDATE RULES

Suppose that $\Theta = \{\theta_1, \ldots, \theta_m\}$ and for all $j$ label $\theta_j$ by $\bar{e}_j \in \mathbb{R}^m$ which is a column vector of all zeros except for its $j$-th element which is equal to one. Furthermore, for all agents $i \in [n]$ let $\mathcal{A}_i$ be the $m$-dimensional probability simplex: $\mathcal{A}_i = \{(x_1, \ldots, x_m)^T \in \mathbb{R}^m : \sum_1^m x_i = 1$ and $x_i \geq 0, \forall i\}$; and for all $\bar{a} := (a_1, \ldots, a_m)^T \in \mathcal{A}_i$ and $\theta_j \in \Theta$, set

$$u_i(\bar{a}, \theta_j) = -\|\bar{a} - \bar{e}_j\|_2^2 := -(1-a_j)^2 - \sum_{\substack{k=1,\\k\neq j}}^m a_k^2.$$

Subsequently, we can calculate the expected pay-off from every such action $\bar{a}$ as:

$$\mathbb{E}_{i,t}\{u_i(\bar{a}, \theta)\} = -\sum_{k=1}^m a_k^2 - 1 + 2\sum_{j=1}^m a_j \boldsymbol{\mu}_{i,t}(\theta_j). \quad (8)$$

Over the $m$-dimensional probability simplex, (8) is uniquely maximized by $\bar{a}^* := (\boldsymbol{\mu}_{i,t}(\theta_1), \ldots, \boldsymbol{\mu}_{i,t}(\theta_m))^T$. Hence, with the probability simplex as their action space $\mathcal{A}_i$ and subjected to the aforementioned utility structure, ever agnet announces her beliefs truthfully, as her optimal action at every epoch of time: $\mathbf{a}_{i,t} = \arg\max_{a \in \mathcal{A}_i} \mathbb{E}_{i,t}\{u_i(a, \theta)\} \equiv \boldsymbol{\mu}_{i,t}(\cdot)$. Therefore, the memoryless update rule $f_i$ in (3) describes how agents' beliefs are being updated following their observations (Fig. 2).

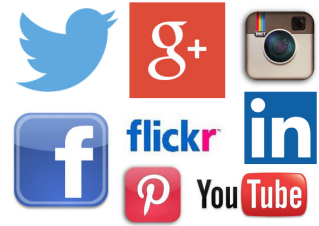

Fig. 2: People reveal their beliefs through status updates and what they share and post on various social media platforms.

In Appendix D, we calculate the following Bayesian belief at time one, in terms of the observed neighboring beliefs and private signal at time zero.

**Lemma 2** (Time-One Bayesian Beliefs). *The Bayesian belief of agent $i$ at time one following her observations of beliefs of her neighbors at time zero and her own private signal at time zero is given by*

$$\boldsymbol{\mu}_{i,1}(\hat{\theta}) = \frac{\nu_i(\hat{\theta}) l_i(\mathbf{s}_{i,0} \mid \hat{\theta}) \left(\prod_{j \in \mathcal{N}(i)} \frac{\boldsymbol{\mu}_{j,0}(\hat{\theta})}{\nu_j(\hat{\theta})}\right)}{\sum_{\tilde{\theta} \in \Theta} \nu_i(\tilde{\theta}) l_i(\mathbf{s}_{i,0} \mid \tilde{\theta}) \left(\prod_{j \in \mathcal{N}(i)} \frac{\boldsymbol{\mu}_{j,0}(\tilde{\theta})}{\nu_j(\tilde{\theta})}\right)}, \forall \hat{\theta} \in \Theta. \quad (9)$$

Subsequently, at any time step $t > 1$, each agent $i$ observes the realized values of $\mathbf{s}_{i,t}$ as well as the current beliefs of her neighbors $\boldsymbol{\mu}_{j,t-1}(\cdot), \forall j \in \mathcal{N}(i)$ and forms a refined opinion $\boldsymbol{\mu}_{i,t}(\cdot)$, using the following rule:

$$\boldsymbol{\mu}_{i,t}(\hat{\theta}) = \frac{\nu_i(\hat{\theta}) l_i(\mathbf{s}_{i,t} \mid \hat{\theta}) \left(\prod_{j \in \mathcal{N}(i)} \frac{\boldsymbol{\mu}_{j,t-1}(\hat{\theta})}{\nu_j(\hat{\theta})}\right)}{\sum_{\tilde{\theta} \in \Theta} \nu_i(\tilde{\theta}) l_i(\mathbf{s}_{i,t} \mid \tilde{\theta}) \left(\prod_{j \in \mathcal{N}(i)} \frac{\boldsymbol{\mu}_{j,t-1}(\tilde{\theta})}{\nu_j(\tilde{\theta})}\right)}, \quad (10)$$

for all $\hat{\theta} \in \Theta$ and at any $t > 1$. In writing (10), every time agent $i$ regards each of her neighbors $j \in \mathcal{N}(i)$ as having started from prior belief $\nu_j(\cdot)$ and arrived at their currently reported belief $\boldsymbol{\mu}_{j,t-1}(\cdot)$ directly, hence rejecting any possibility of a past history. This is equivalent to the assumption that the reported beliefs of every neighbor are formed from a private observation and a fixed prior, and not through repeated communications.

Such a rule is of course not the optimum Bayesian update of agent $i$'s belief at any step $t > 1$, because the agent is not taking into account the complete observed history of her private signals and neighbors' beliefs and is instead, basing her inference entirely on the immediately observed signal and neighboring beliefs; hence, the name *memoryless*. Here, the status of a *Rational but Memoryless* agent is akin to a person who is possessed of a knowledge but cannot see how she has come to be possessed of that knowledge. Likewise, it is by the requirement of rationality in such a predicament that we impose a fixed prior $\nu_i(\cdot)$ on every agent $i$ and carry it through for all times $t$. Indeed, it is the grand tradition of Bayesian statistics, as advocated in the prominent and influential works of [76], [77], [78], [79] and many others, to argue on normative grounds that rational behavior in a decision theoretic framework forces individuals to employ Bayes rule and appropriate it to their personal priors.

*A. Analysis of Convergence and Log-Linear Learning*

A main question of interest is whether the agents can learn the true realized value $\theta$:

**Definition 1** (Learning). *An agent $i$ is said to learn the truth, if $\lim_{t\to\infty} \boldsymbol{\mu}_{i,t}(\theta) = 1$, $\mathbb{P}_\theta$-almost surely.*

We begin our analysis of convergence and learning under the update rule in (10) by considering the case of a single agent $i$, who starts from a prior belief $\nu_i(\cdot)$ and sequentially updates her beliefs according to Bayes rule:

$$\boldsymbol{\mu}_{i,t}(\hat{\theta}) = \frac{\boldsymbol{\mu}_{i,t-1}(\hat{\theta})\ell_i(\mathbf{s}_{i,t} \mid \hat{\theta})}{\sum_{\tilde{\theta}\in\Theta}\boldsymbol{\mu}_{i,t-1}(\tilde{\theta})\ell_i(\mathbf{s}_{i,t} \mid \tilde{\theta})}, \forall \hat{\theta} \in \Theta. \quad (11)$$

The Bayesian belief update in (11) linearizes in terms of the log-ratio of beliefs and signal likelihoods, $\phi_{i,t}(\cdot)$ and $\lambda_{\check{\theta}}(\cdot)$, leading to

$$\phi_{i,t}(\check{\theta}) = \log\left(\frac{\nu_i(\check{\theta})}{\nu_i(\theta)}\right) + \sum_{\tau=0}^{t}\lambda_{\check{\theta}}(\mathbf{s}_{i,\tau})$$
$$\to \log\left(\frac{\nu_i(\check{\theta})}{\nu_i(\theta)}\right) + (t+1)\mathbb{E}_\theta\{\lambda_{\check{\theta}}(\mathbf{s}_{i,0})\} \quad (12)$$

$\mathbb{P}_\theta$-almost surely, as $t \to \infty$; by the strong law of large numbers [80, Theorem 22.1] applied to the sequence of $\mathbb{E}_\theta$-integrable, independent and identically distributed variables: $\lambda_{\check{\theta}}(\mathbf{s}_{i,t}), t \in \mathbb{N}_0$. In particular, if $D_{KL}\left(\ell_i(\cdot|\theta)||\ell_i(\cdot|\check{\theta})\right) := -\mathbb{E}_\theta\{\lambda_{\check{\theta}}(\mathbf{s}_{i,t})\} > 0$, then $\phi_{i,t}(\check{\theta}) \to -\infty$ almost surely and agent $i$ asymptotically rejects the false state $\check{\theta}$ in favor of the true state $\theta$, putting a vanishing belief on the former relative to the latter. Therefore, the single Bayesian agent following (11) learns the truth if and only if $D_{KL}\left(\ell_i(\cdot|\theta)||\ell_i(\cdot|\check{\theta})\right) > 0$ for all $\check{\theta} \neq \theta$ and the learning is asymptotically exponentially fast at the rate $\min_{\check{\theta}\in\Theta\setminus\{\theta\}} D_{KL}\left(\ell_i(\cdot|\theta)||\ell_i(\cdot|\check{\theta})\right)$ as shown in [81].[1]

The preceding result is also applicable to the case of a Bayesian agents with direct (centralized) access to all observations across the network: consider an outside Bayesian agent $\hat{o}$ who shares the same common knowledge of the prior and signal structures with the networked agents; in particular, $\hat{o}$ knows the signal structures $\ell_i(\cdot|\hat{\theta})$, for all $\hat{\theta} \in \Theta$ and $i \in [n]$; thence, making the same inferences as any other agent when given access to the same observations. Consider next a Gedanken experiment where $\hat{o}$ is granted direct access to all the signals of every agent at all times. The analysis leading to (12) can be applied to the evolution of $\log$ belief ratios for $\hat{o}$, whose observations at every time $t \in \mathbb{N}_0$ is an element of the product space $\prod_{i\in[n]}\mathcal{S}_i$. Subsequently, the centralized Bayesian beliefs concentrate on the true state at the asymptotically exponentially fast rate of

$$R_n := \min_{\check{\theta}\in\Theta\setminus\{\theta\}} D_{KL}\left(\prod_{i\in[n]}\ell_i(\cdot|\theta) \Big\| \prod_{i\in[n]}\ell_i(\cdot|\check{\theta})\right)$$
$$= \min_{\check{\theta}\in\Theta\setminus\{\theta\}} \sum_{i\in[n]} D_{KL}\left(\ell_i(\cdot|\theta)||\ell_i(\cdot|\check{\theta})\right). \quad (13)$$

Next to understand the evolution of beliefs under the log-linear updates in (10), consider the network graph structure as encoded by its adjacency matrix $A$ defined as $[A]_{ij} = 1 \iff (j,i) \in \mathcal{E}$, and $[A]_{ij} = 0$ otherwise. For a strongly connected $\mathcal{G}$ the Perron-Frobenius theory [83, Theorem 1.5] implies that $A$ has a simple positive real eigenvalue, denoted by $\rho > 0$, which is equal to its spectral radius. Moreover, the left eigenspace associated with $\rho$ is one-dimensional with the corresponding eigenvector $\overline{\alpha} = (\alpha_1,\ldots,\alpha_n)^T$, uniquely satisfying $\sum_{i=1}^n \alpha_i = 1$, $\alpha_i > 0$, $\forall i \in [n]$, and $\overline{\alpha}^T A = \rho\overline{\alpha}^T$. The entry $\alpha_i$ is also called the centrality of agent $i$ and as the name suggests, it is a measure of how central is the location of agent in the network. Our main result state that almost sure learning cannot be realized in a strongly connected network unless it has unit spectral radius which is the case only of a directed circle.

**Theorem 1** (No Learning when Spectral Radius $\rho > 1$). *In a strongly connected social network and under the memoryless belief updates in* (10)*, no agents can learn the truth unless the*

---

[1]Note from the information inequality for the Kullback-Leibler divergence that $D_{KL}(\cdot||\cdot) \geq 0$ and the inequality is strict whenever $\ell_i(\cdot|\check{\theta}) \not\equiv \ell_i(\cdot|\theta)$, i.e. $\exists s \in \mathcal{S}_i$ such that $\ell_i(s|\check{\theta}) \neq \ell_i(s|\theta)$ [82, Theorem 2.6.3]. Further note that whenever $\ell_i(\cdot|\check{\theta}) \equiv \ell_i(\cdot|\theta)$ or equivalently $D_{KL}\left(\ell_i(\cdot|\theta)||\ell_i(\cdot|\check{\theta})\right) = 0$, then the two states $\check{\theta}$ and $\theta$ are statically indistinguishable to agent $i$: there is no way for agent $i$ to distinguish between $\check{\theta}$ and $\theta$, based only on her received signals. This is because both $\theta$ and $\check{\theta}$ induce the same probability distribution on her sequence of observed i.i.d. signals. Since different states $\hat{\theta} \in \Theta$ are distinguished through their different likelihood functions $\ell_i(\cdot \mid \hat{\theta})$; the more refined such differences are, the better the states are distinguished. Hence, the proposed asymptotic rate is one measure of resolution for the likelihood structure of agent $i$.

*spectral radius* $\rho = 1$.

**Proof outline:** A complete proof is included in Appendix E, but here we provide a description of the mechanism and the interplay between the belief aggregation and information propagation. To facilitate the exposition of the underlying logic we introduce some notation. We define a global (network-wide) random variable $\boldsymbol{\Phi}_t(\check{\theta}) := \sum_{i=1}^{n} \alpha_i \phi_{i,t}(\check{\theta})$, where $\alpha_i$ is the centrality of agent $i$ and $\boldsymbol{\Phi}_t(\check{\theta})$ characterizes how biased (away from the truth and towards $\check{\theta}$) the network beliefs and priors are at each point in time. In particular, if any agent is to learn the truth, then $\boldsymbol{\Phi}_t(\check{\theta}) \to -\infty$ as $t \to \infty$ for all the false states $\check{\theta} \in \Theta \setminus \{\theta\}$. To proceed, we define another network-wide random variable $\boldsymbol{\Lambda}_t(\check{\theta}) := \sum_{i=1}^{n} \alpha_i \lambda_{\check{\theta}}(\mathbf{s}_{i,t})$ which characterizes the information content of the observed signals (received information) for the entire network, at each time $t$. Moreover, since the received signal vectors $\{(\mathbf{s}_{1,t}, \ldots, \mathbf{s}_{n,t}), t \in \mathbb{N}_0\}$ are i.i.d. over time, $\forall \check{\theta} \neq \theta$, $\{\boldsymbol{\Lambda}_t(\check{\theta}), t \in \mathbb{N}_0\}$ constitutes a sequence of i.i.d. random variables satisfying $\mathbb{E}\{\boldsymbol{\Lambda}_t(\check{\theta})\} = -\sum_{i=1}^{n} \alpha_i D_{KL}\left(\ell_i(\cdot|\theta) || \ell_i(\cdot|\check{\theta})\right) \leqslant 0$. In order for the agents to learn the true state of the world based on their observations, it is necessary that at each false state $\check{\theta} \neq \theta$ some agent be able to distinguish $\check{\theta}$ from the truth $\theta$, in which case $\mathbb{E}\{\boldsymbol{\Lambda}_t(\check{\theta})\} < 0$, and we can refer to this criterion as *global identifiablity* for the true state $\theta$.[1]

In Appendix E we argue that under the update rules in (10) the global belief ratio statistics $\boldsymbol{\Phi}_t(\check{\theta})$ evolves as a sum of weighted i.i.d. variables $\rho^\tau \boldsymbol{\Lambda}_{t-\tau}(\check{\theta})$:

$$\boldsymbol{\Phi}_t(\check{\theta}) = \sum_{\tau=0}^{t} \rho^\tau \left(\boldsymbol{\Lambda}_{t-\tau}(\check{\theta}) + (1-\rho)\beta(\check{\theta})\right), \qquad (14)$$

where $\beta(\check{\theta}) := \sum_{i=1}^{n} \alpha_i \log\left(\nu_i(\check{\theta})/\nu_i(\theta)\right)$ is a measure of bias in the initial prior beliefs. The weights in (14) form a geometric progressions in $\rho$; hence, the variables increase unbounded in their variance and convergence cannot hold true in a strongly connected social network, unless $\rho = 1$. This is due to the fact that $\rho$ upper bounds the average degree of the graph [84, Chapter 2], and every node in a strongly connected graph has degree greater than or equal to one, subsequently $\rho \geq 1$ for all strongly connected graphs. □

**Remark 2** (Polarization, data incest and unlearning). *The unlearning in the case of $\rho > 1$ in Theorem 1, which applies to all strongly connected topologies except directed circles (where $\rho = 1$, see Subsection IV-B below), is related to the inefficiencies associated with social learning and can be attributed to the agents' naivety in inferring the sources of their information, and their inability to interpret the actions of their neighbors rationally [85]. In particular, when $\rho > 1$ the noise or randomness in the agents' observations is amplified at every stage of network interactions; since the agents fail to correct for the repetitions in the sources of their observations as in the case of persuasion bias argued by DeMarzo, Vayanos and Zwiebel [9], or data incest argued by Krishnamurthy and Hoiles [11]. When $\rho > 1$ the effect of the agents' priors is also amplified through the network interactions and those states $\hat{\theta}$ for which $\beta(\hat{\theta}) > 0$ in (14), will be asymptotically rejected as $\sum_{\tau=0}^{t} \rho^\tau(1-\rho)\beta(\check{\theta}) \to -\infty$, irrespectively of the observed data $\boldsymbol{\Lambda}_\tau(\check{\theta}), \tau \in \mathbb{N}_0$. This phenomenon arises as agents engage in excessive anti-imitative behavior, compensating for the neighboring priors at every period [44]. It is justified as a case of choice shift toward more extreme opinions [5], [6] or group polarization [7], [8], when like-minded people after interacting with each other and under the influence of their mutually positive feedback become more extreme in their opinions, and less receptive of opposing beliefs.*

### B. Learning in Circles and General Connected Topologies

For a strongly connected digraph $\mathcal{G}$, if $\rho = 1$, then it has to be the case that all nodes have degree 1 and the graph is a directed circle. Subsequently, the progression for $\boldsymbol{\Phi}_t(\check{\theta})$ in (14) reduces to sum of i.i.d. variables in $\mathcal{L}^1$ and by the strong law of large numbers [80, Theorem 22.1], it converges almost surely to the mean value

$$\boldsymbol{\Phi}_t(\check{\theta}) = \beta(\check{\theta}) + \sum_{\tau=0}^{t} \boldsymbol{\Lambda}_\tau(\check{\theta}) \to \beta(\check{\theta}) + (t+1)\mathbb{E}\{\boldsymbol{\Lambda}_0(\check{\theta})\}$$
$$\to -\infty,$$

as $t \to \infty$, provided that $\mathbb{E}\{\boldsymbol{\Lambda}_0(\check{\theta})\} < 0$, i.e. if the truth is globally identifiable. Note also the analogy with (12), where $\boldsymbol{\Lambda}_t(\check{\theta})$ is replaced by $\lambda_{\check{\theta}}(\mathbf{s}_{i,t})$ as both represent the observed signal(s) or received information at time $t$. Indeed, if we further assume that $\boldsymbol{\nu}_i(\cdot) \equiv \nu(\cdot)$ for all $i$, i.e. all agents share the same common prior, then (10) for a circular network becomes

$$\boldsymbol{\mu}_{i,t}(\hat{\theta}) = \frac{\boldsymbol{\mu}_{j,t-1}(\hat{\theta})\ell_i(\mathbf{s}_{i,t} \mid \hat{\theta})}{\sum_{\tilde{\theta} \in \Theta} \boldsymbol{\mu}_{j,t-1}(\tilde{\theta})\ell_i(\mathbf{s}_{i,t} \mid \tilde{\theta})}, \forall \hat{\theta} \in \Theta, \qquad (15)$$

where $j \in [n]$ is the unique vertex $j \in \mathcal{N}(i)$. Update (15) replicates the Bayesian update of a single agents in (11) but the self belief $\boldsymbol{\mu}_{i,t-1}(\cdot)$ on the right-hand side being is replaced by the belief $\boldsymbol{\mu}_{j,t-1}(\cdot)$ of the unique neighbor $\{j\} = \mathcal{N}(i)$. Indeed, the learning in this case is asymptotically exponentially fast at the rate $(1/n) \min_{\check{\theta} \in \Theta \setminus \{\theta\}} \sum_{j=1}^{n} D_{KL}\left(\ell_j(\cdot|\theta)||\ell_j(\cdot|\check{\theta})\right) = 1/3R_3$; hence, the same exponential rate as that of a central Bayesian can be achieved through the BWR update rule, except for a $1/n$ factor that decreases with the increasing cycle length, cf. [81].

**Example 1** (Eight Agents with Binary Signals in a Tri-State World.). *Consider the network of agents in Fig. 3 with the true state of the world being 1, the first of the tree possible states $\Theta = \{1, 2, 3\}$. The agents receive binary signals about the true state $\theta$ according to the likelihoods listed in the table.*

---

[1] The global identifiability condition can be also viewed in the following sense: consider a gedanken experiment where an external fully rational observer $\hat{o}$ is granted direct access to all the signals of all agents in the network and assume further that she shares the same common knowledge of the prior and signal structures with the network agents. Then $\hat{o}$ learns the truth if, and only if, it is globally identifiable.

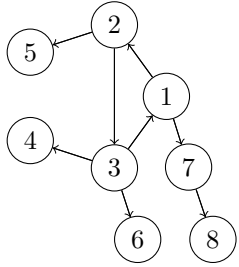

| likelihoods | $\hat{\theta}=1$ | $\hat{\theta}=2$ | $\hat{\theta}=3$ |
|---|---|---|---|
| $l_1(\mathbf{s}_{1,t}=0 \mid \hat{\theta})$ | $\frac{1}{3}$ | $\frac{1}{3}$ | $\frac{1}{5}$ |
| $l_2(\mathbf{s}_{2,t}=0 \mid \hat{\theta})$ | $\frac{1}{2}$ | $\frac{2}{3}$ | $\frac{1}{2}$ |
| $l_3(\mathbf{s}_{3,t}=0 \mid \hat{\theta})$ | $\frac{1}{4}$ | $\frac{1}{4}$ | $\frac{1}{4}$ |

Fig. 3: A hybrid structure

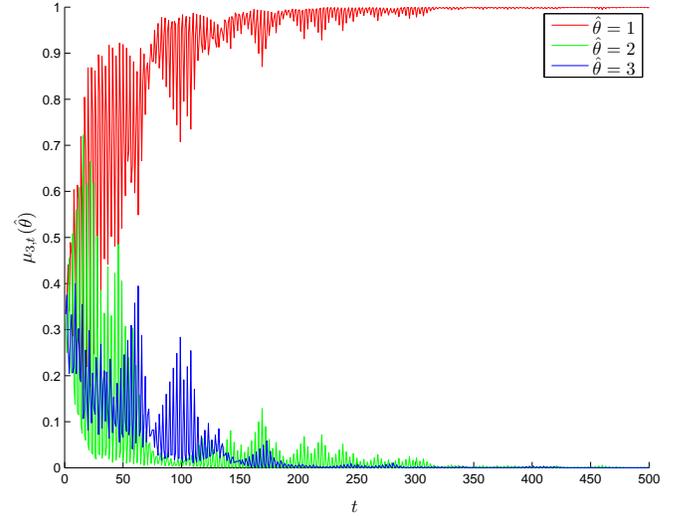

Fig. 4: Evolution of the third agent's beliefs over time

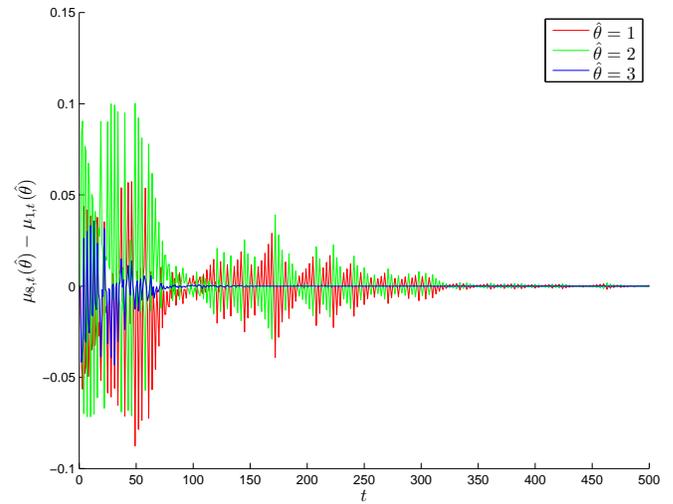

Fig. 5: The difference between the first and eighth agents' beliefs over time

*We begin by the observation that this network can be thought of as a rooted directed tree, in which the root node is replaced with a directed circle (the root circle).[1] Next note that the root circle is comprised of three agents and none of them can learn the truth on their own. Indeed, agent 3 does not receive any informative signals; therefore, in isolation i.e. using (11), her beliefs shall never depart from their initial priors. We further set $l_j(\cdot \mid \cdot) \equiv l_3(\cdot \mid \cdot)$ for all $j \in [8]\setminus[3]$, so that all the peripheral follower agents are also unable to infer anything about the true state of the world from their own private signals.*

*Starting from a uniform common prior and following the proposed rules (15), all agents asymptotically learn the true state, even though none of them can learn the true state on their own. The plots in Figs. 4 and 5 depict the evolutions of the beliefs for the third agent as well as the difference between the beliefs for the first and eighth agents. We can further show that all agents learn the true state at the same exponentially fast asymptotic rate. In fact, the three nodes belonging to the directed circle learn the true state of the world at the exponentially fast asymptotic rate of $(1/3)R_3$ noted above, irrespectively of the peripheral nodes. The remaining peripheral nodes then follow up with the beliefs of root circle nodes, except for a vanishing difference that increases with the increasing distance of a peripheral node from the root circle: following (15), the first three agents form a circle of leaders where they combine their observations and reach a consensus; every other agent in the network then follows whatever state that the leaders have collectively agreed upon.*

In [86] the authors show the application of the update rule in (15) to general strongly connected topologies where agents have more than just a single neighbor in their neighborhoods. It is proposed to choose a neighbor $j \in \mathcal{N}(i)$ independently at random every time and then apply (15) with the reported belief from that neighbor. Here again if the truth is globally identifiable, all agents learn the truth at an asymptotically exponentially fast rate given by $\min_{\check{\theta} \in \Theta \setminus \{\theta\}} \sum_{j=1}^{n} \pi_j D_{KL}\left(\ell_j(\cdot|\theta)||\ell_j(\cdot|\check{\theta})\right)$, where $\pi_j$ are the probabilities in the stationary distribution of the Markov chain whose transition probabilities are the same as the probabilities for the random choice of neighbors at every point in time.[1] It is notable that the asymptotic rate here is a weighted average of the $KL$ distances $D_{KL}\left(\ell_i(\cdot|\theta)||\ell_i(\cdot|\check{\theta})\right)$, in contrast with the arithmetic (unweighted) mean $(1/nR_n)$ that arise in the circular case. Both rates are upper bounded by the centralized Bayesian learning rate of $R_n$ calculated in (13). Finally, we point out that the rate of distributed learning upper bounds the (weighted) average of individual learning rates. It is due to the fact that observations of different agents complement each other, and while one agent may be good at distinguishing one false state from the truth, she can rely on observational abilities of other agents for distinguishing the remaining false states: consider agents 1 and 2 in Example 1, the former can distinguish $\hat{\theta} = 2$ from $\theta = 1$, while the latter is good at distinguishing $\hat{\theta} = 3$ from $\theta = 1$; together they can distinguish all states. Hence, the overall rate of distributed learning upper bounds the average of individual learning rates, and is itself upper bounded by the learning rate of a central Bayesian agent:

$$\frac{1}{n}\sum_{i=1}^{n} \min_{\check{\theta}\in\Theta\setminus\{\theta\}} D_{KL}\left(\ell_i(\cdot|\theta)||\ell_i(\cdot|\check{\theta})\right)$$
$$< \min_{\check{\theta}\in\Theta\setminus\{\theta\}} \frac{1}{n}\sum_{i=1}^{n} D_{KL}\left(\ell_i(\cdot|\theta)||\ell_i(\cdot|\check{\theta})\right) = \frac{1}{n}R_n < R_n.$$

---

[1] Any weakly connected digraph $\mathcal{G}$ which has only degree zero or degree one nodes can be drawn as a rooted tree, whose root is replaced by a directed circle, a so-called root circle. This is true since any such digraph can have at most one directed circle and all other nodes that are connected to this circle should be directed away from it, otherwise $\mathcal{G}$ would have to include a node of degree two or higher.

Fixing the priors over time will not result in convergence of beliefs, except in very specific cases as discussed above. In [90] we investigate the properties of convergence and learning under the update rules in (10), where the priors $\nu_j(\cdot)$ are replaced by time-varying distributions $\boldsymbol{\xi}_{i,j}(\cdot,t)$ that parametrize the log-linear updating of the agents' beliefs over time. It is notable that the memoryless Bayesian update in (10) has a log-linear structure similar to non-Bayesian update rules studied in the literature [91]–[95]; the roots for such a geometric averaging of the neighboring beliefs is traced to logarithmic opinion pools [96], [97] and can be also justified under specific behavioral assumptions [53].

## V. CONCLUSIONS

This work addressed a social and observational learning model in social networks. Agents attempt to learn some unknown state of the world which belongs to a finite state space. Conditioned on the true state, a sequence of i.i.d. private signals are generated and observed by each agent of the network. The private signals do not provide each agent with adequate information to identify the truth. Hence, agents interact with their neighbors to augment their imperfect observations with those of their neighbors. We proposed a belief aggregation and inference scheme that we call Bayesian without Recall (BWR), as a behavioral model to interpret and provide justification for the variety of non-Bayesian update rules that are suggested in the literature. Accordingly, by replicating the rule that maps the initial priors, neighbor's decisions, and the private signal to Bayesian posterior at one time step for all future time steps, one can derive non-Bayesian updates with well-known and intuitive structures, such as majority rules or log-linear belief updates.

Following the BWR approach, the complexities of a fully rational inference at the forthcoming epochs are avoided, while some essential features of Bayesian inference are preserved. We analyzed the specific form of BWR updates in two cases of binary state and action space, as well as a finite state space with actions taken over the probability simplex. In the case of binary actions the BWR updates take the form of a linear majority rule, whereas if the action spaces are rich enough for the agents to reveal their beliefs, then belief updates take a log-linear format. In each case we investigate the properties of convergence, consensus and learning; the latter is of particular interest, in a strongly connected social network when the truth is identifiable through the aggregate private observations of all individuals but not individually.

On the one hand, the specific forms of the BWR update rules in each case help us to better understand the mechanisms of naive inference, when rational agents are devoid of their ability to make recollections. On the other hand, our results also highlight the consequences of such naivety in shaping the mass behavior; by comparing our predications with the rational learning outcomes. In particular, we saw in Subsection III-B that there is a positive probability for rational but memoryless agents in an Ising model to mis-learn by reaching consensus on an untruth. However Bayesian (fully rational) beliefs constitute a bounded martingale; hence, when truth is identifiable and number of observations increases, the beliefs of rational agents converge almost surely to a point mass centered at the true state [28], [29]. Similarly Theorem 1 states the impossibility of asymptotic learning under the BWR belief updates, whenever the spectral radius of the interconnection graph adjacency is greater than one.

Last but not least, we pinpoint a key difference between the BWR action and belief updates: the former are weighted updates, wheras the latter are unweighted symmetric updates. Accordingly, an agent weighs each neighbor's action differently and in accordance with the quality of private signals which are inferred from actions. On the other hand, when communicating their beliefs the quality of each neighbor's signal is already internalized in their reported belief; hence, when incorporating her neighboring beliefs, an agent regards the reported beliefs of all her neighbors equally, and irrespective of the quality of their private signals.

## APPENDIX A
## PROOF OF LEMMA 1: TIME-ONE BAYESIAN ACTIONS

Note that given her observation of private signal $\mathbf{s}_{i,0}$, the posterior probability assigned by agent $i$ to the set $\theta_1$ is given by (1) with $\hat{\theta} = \theta_1$. We form a dichotomy of the signal space $\mathcal{S}_i$ of each agent into $\mathcal{S}_i^1$ and $\mathcal{S}_i^{-1}$; by setting $\mathcal{S}_i^1 := \{s \in \mathcal{S}_i : \ell_i(s|\theta_1)\nu_i(\theta_1) \geq \ell_i(s|\theta_2)\nu_i(\theta_2)\}$ and $\mathcal{S}_i^{-1} := \mathcal{S}_i \setminus \mathcal{S}_i^1$. It thus follows from (6) that for any $j \in \mathcal{N}(i)$ the observation that $\mathbf{a}_{j,0} = 1$ is equivalent to the information that $\{\mathbf{s}_{j,0} \in \mathcal{S}_j^1\}$ and $\mathbf{a}_{j,0} = -1$ is equivalent to the information that $\{\mathbf{s}_{j,0} \in \mathcal{S}_j^{-1}\}$. Thereby, the belief of agent $i$ at time $t = 1$ given her observation of the actions of her neighbors and the private signal $\mathbf{s}_{j,0}$ is given by

$$\boldsymbol{\mu}_{i,1}(\theta_1) = \frac{\ell_i(\mathbf{s}_{i,0}|\theta_1) \prod_{j \in \mathcal{N}(i)} \left( \sum_{s_j \in \mathcal{S}_j^{\mathbf{a}_{j,0}}} \ell_j(s_j|\theta_1) \right) \nu_i(\theta_1)}{\sum_{\hat{\theta} \in \Theta} \ell_i(\mathbf{s}_{i,0}|\hat{\theta}) \prod_{j \in \mathcal{N}(i)} \left( \sum_{s_j \in \mathcal{S}_j^{\mathbf{a}_{j,0}}} \ell_j(s_j|\hat{\theta}) \right) \nu_i(\hat{\theta})},$$

and we can thus form the ratio

$$\frac{\boldsymbol{\mu}_{i,1}(\theta_1)}{\boldsymbol{\mu}_{i,1}(\theta_2)} = \frac{\ell_i(\mathbf{s}_{i,0}|\theta_1)\nu_i(\theta_1)}{\ell_i(\mathbf{s}_{i,0}|\theta_2)\nu_i(\theta_2)} \prod_{j \in \mathcal{N}(i)} \left( \frac{\sum_{s_j \in \mathcal{S}_j^{\mathbf{a}_{j,0}}} \ell_j(s_j|\theta_1)}{\sum_{s_j \in \mathcal{S}_j^{\mathbf{a}_{j,0}}} \ell_j(s_j|\theta_2)} \right)$$

$$= \frac{\ell_i(\mathbf{s}_{i,0}|\theta_1)\nu_i(\theta_1)}{\ell_i(\mathbf{s}_{i,0}|\theta_2)\nu_i(\theta_2)} V_i \prod_{j \in \mathcal{N}(i)} W_j^{\mathbf{a}_{j,0}}, \quad (16)$$

where for all $i \in [n]$ we have defined

$$V_i = \prod_{j \in \mathcal{N}(i)} \left( \frac{\sum_{s_j \in \mathcal{S}_j^1} \ell_j(s_j|\theta_1)}{\sum_{s_j \in \mathcal{S}_j^1} \ell_j(s_j|\theta_2)} \times \frac{\sum_{s_j \in \mathcal{S}_j^{-1}} \ell_j(s_j|\theta_1)}{\sum_{s_j \in \mathcal{S}_j^{-1}} \ell_j(s_j|\theta_2)} \right)^{1/2} \quad (17)$$

$$W_i = \left( \frac{\sum_{s_i \in \mathcal{S}_i^1} \ell_i(s_i|\theta_1)}{\sum_{s_i \in \mathcal{S}_i^1} \ell_i(s_i|\theta_2)} \times \frac{\sum_{s_i \in \mathcal{S}_i^{-1}} \ell_i(s_i|\theta_2)}{\sum_{s_i \in \mathcal{S}_i^{-1}} \ell_i(s_i|\theta_1)} \right)^{1/2}. \quad (18)$$

---

[1] In many distributed learning models over random and switching networks, agents must have positive self-reliant at any time; as for instance in gossip algorithms [87] and ergodic stationary processes [88]. This condition however is relaxed under (15), as our agents rely entirely on the beliefs of their neighbors every time that they select a neighbor to communicate with. Moreover, unlike the majority of results that rely on the convergence properties of products of stochastic matrices and are applicable only to irreducible and aperiodic communication matrices, cf. [10, Proposition 1]; the convergence results in [86] do not require the transition probability matrix to be aperiodic, as it relies on properties of ergodic Markov chains and holds true for any irreducible, finite-state chain [89, Theorems 1.5.6 and 1.7.7].

Furthermore let $w_i := \log W_i$ and $\eta_i := \log(V_i \nu_i(\theta_1)/\nu_i(\theta_2))$ be constants that are determined completely by the initial prior and signal structures of each agent and her neighbors. Subsequently, taking logarithms of both sides in (16) yields the following update rule for the log-ratio of the beliefs at time one,

$$\log\left(\frac{\boldsymbol{\mu}_{i,1}(\theta_1)}{\boldsymbol{\mu}_{i,1}(\theta_2)}\right) = \sum_{j \in \mathcal{N}(i)} w_j \mathbf{a}_{j,0} + \eta_i + \lambda_1(\mathbf{s}_{i,0}). \quad (19)$$

Finally, we can apply (5) to derive the claimed expression in Lemma 1 for the updated Bayesian action of agent $i$ following her observations of her neighbors' actions $\mathbf{a}_{j,0}$, $j \in \mathcal{N}(i)$ and her own private signal $\mathbf{s}_{i,0}$. We end our derivation by pointing out some facts concerning constants $\eta_i$ and $w_i$ which appear in (19).

**Lemma 3** (Monotonicity of $\eta_i$). *Consider any $i \in [n]$ and fix a signal $s_j \in \mathcal{S}_j$ for some $j \in \mathcal{N}(i)$. It holds true that the constant $\eta_i$ is increasing in $\ell_j(s_j \mid \theta_1)$ and decreasing in $\ell_j(s_j \mid \theta_2)$.*

*Proof.* The claim follows directly from the defining relation $\eta_i = \log(\nu_i(\theta_1)/\nu_i(\theta_2)) + \log V_i$, as replacing from (17) yields

$$\log V_i = +\frac{1}{2} \sum_{j \in \mathcal{N}(i)} \log \sum_{s_j \in \mathcal{S}_j^1} \ell_j(s_j|\theta_1)$$

$$+ \frac{1}{2} \sum_{j \in \mathcal{N}(i)} \log \sum_{s_j \in \mathcal{S}_j^{-1}} \ell_j(s_j|\theta_1)$$

$$- \frac{1}{2} \sum_{j \in \mathcal{N}(i)} \log \sum_{s_j \in \mathcal{S}_j^1} \ell_j(s_j|\theta_2)$$

$$- \frac{1}{2} \sum_{j \in \mathcal{N}(i)} \log \sum_{s_j \in \mathcal{S}_j^{-1}} \ell_j(s_j|\theta_2). \quad (20)$$

The proof now follows upon the realization that for any fixed $s_j \in \mathcal{S}_j$, $j \in \mathcal{N}(i)$ the term $\ell_j(s_j \mid \theta_1)$ appears in one of the first two terms appearing with a plus sign in (20), and the term $\ell_j(s_j \mid \theta_2)$ appears in one of the last first two terms appearing with a minus sign in (20). Hence, when all else kept constant, $\log V_i$ and subsequently $\eta_i$ is increasing in $\ell_j(s_j \mid \theta_1)$ and decreasing in $\ell_j(s_j \mid \theta_2)$. □

**Lemma 4** (Positivity of $w_i$). *It holds true for any $i \in [n]$ that $w_i \geq 0$.*

*Proof.* First note from the definitions of the sets $\mathcal{S}_i^1$ and $\mathcal{S}_i^{-1}$ that $\forall s \in \mathcal{S}_i^1$,

$$\frac{\ell_i(s|\theta_1)}{\ell_i(s|\theta_2)} \geq \frac{\nu_i(\theta_2)}{\nu_i(\theta_1)}, \text{ and } \forall s \in \mathcal{S}_i^{-1}, \frac{\ell_i(s|\theta_2)}{\ell_i(s|\theta_1)} > \frac{\nu_i(\theta_1)}{\nu_i(\theta_2)}.$$

Next we sum the numerators and denominators of the likelihood ratios of the signals in each of sets $\mathcal{S}_i^1$ and $\mathcal{S}_i^{-1}$; invoking basic algebraic properties from the resultant fractions yields

$$\frac{\sum_{s \in \mathcal{S}_i^1} \ell_i(s|\theta_1)}{\sum_{s \in \mathcal{S}_i^1} \ell_i(s|\theta_2)} \geq \frac{\nu_i(\theta_2)}{\nu_i(\theta_1)}, \text{ and } \frac{\sum_{s \in \mathcal{S}_i^{-1}} \ell_i(s|\theta_2)}{\sum_{s \in \mathcal{S}_i^{-1}} \ell_i(s|\theta_1)} > \frac{\nu_i(\theta_1)}{\nu_i(\theta_2)}.$$

Subsequently, replacing form (18) yields that

$$W_i := \frac{\sum_{s \in \mathcal{S}_i^1} \ell_i(s|\theta_1)}{\sum_{s \in \mathcal{S}_i^1} \ell_i(s|\theta_2)} \times \frac{\sum_{s \in \mathcal{S}_i^{-1}} \ell_i(s|\theta_2)}{\sum_{s \in \mathcal{S}_i^{-1}} \ell_i(s|\theta_1)} \geq \frac{\nu_i(\theta_2)}{\nu_i(\theta_1)} \times \frac{\nu_i(\theta_1)}{\nu_i(\theta_2)}$$

$$= 1,$$

and proof follows from the defining relation $w_i := \log W_i \geq 0$. □

## APPENDIX B
## A MARKOV CHAIN ON THE BOOLEAN CUBE

To begin, for any vertex of the Boolean hypercube $\bar{a} := (a_1, \ldots, a_n)^T \in \{\pm 1\}^n$ and each agent $i$, define the function $\pi_i : \{\pm 1\}^n \to [0, 1]$ as $\pi_i(\bar{a}) := \mathbb{P}\{\bar{\mathbf{a}}_{i,t+1} = +1 \mid \bar{\mathbf{a}}_t = $

$\overline{a}\} = \mathbb{P}_\theta\{-\lambda_1(\mathbf{s}_{i,t+1}) \leq \sum_{j\in\mathcal{N}(i)} w_j a_j + \eta_i\}$. The transition probabilities for the Markov chain of action profiles on the Boolean hypercube are given by

$$P(\overline{a}', \overline{a}) := \mathbb{P}\{\overline{\mathbf{a}}_{t+1} = \overline{a}' \mid \overline{\mathbf{a}}_t = \overline{a}\}$$
$$= \prod_{i:a'_i=+1} \pi_i(\overline{a}) \prod_{i:a'_i=-1} (1 - \pi_i(\overline{a})), \quad (21)$$

for all $t \in \mathbb{N}_0$ and any pair of vertices $\overline{a}' := (a'_1, \ldots, a'_n)^T \in \{\pm 1\}^n$ and $\overline{a} \in \{\pm 1\}^n$.

It follows from the classification of states and chains in [98, Section 2.4] that $\{\pm 1\}^n$ can be partitioned into sets of transient communication classes: $\mathcal{C}'_1, \ldots, \mathcal{C}'_{r'}$, and recurrent (ergodic) communication classes: $\mathcal{C}_1, \ldots, \mathcal{C}_r$. Moreover, as $t \to \infty$, $\overline{\mathbf{a}}_t$ almost surely belongs to $\cup_{i\in[r]}\mathcal{C}_i$. It is further true that if $\overline{\mathbf{a}}_{t_0} \in \mathcal{C}_i$ for some $i \in [r]$ and $t_0 \in \mathbb{N}$, then $\overline{\mathbf{a}}_t \in \mathcal{C}_i$ almost surely for all $t \geq t_0$: the process will almost surely leave any set of transient action profiles, i.e. $\cup_{i\in[r']}\mathcal{C}'_i$, and will almost surely remain in the first recurrent set that it reaches before any other. Let $\mathbf{r}^* := \arg\min_{\rho\in[r]}\{t : \overline{\mathbf{a}}_t \in \mathcal{C}_\rho\}$ be the random variable that determines the first ergodic set of action profiles that is reached by the Markov chain process $\{\overline{\mathbf{a}}_t, t \in \mathbb{N}_0\}$; suppose $\boldsymbol{\tau} := \text{card}(\mathcal{C}_{\mathbf{r}^*})$ and further denote $\mathcal{C}_{\mathbf{r}^*} := \{\mathbf{a}^*_1, \ldots, \mathbf{a}^*_{\boldsymbol{\tau}}\}$. The asymptotic behavior of the process can now be characterized as follows.

**Proposition 2** (Asymptotic Distribution of Action Profiles). *Let $\overline{\mathbf{p}} := (\mathbf{p}_1, \ldots, \mathbf{p}_{\boldsymbol{\tau}})^T$ be the stationary distribution over $\mathcal{C}_{\mathbf{r}^*}$ which uniquely satisfies $\mathbf{p}_k \sum_{j=1}^{\boldsymbol{\tau}} P(\mathbf{a}^*_k, \mathbf{a}^*_j) = \sum_{j=1}^{\boldsymbol{\tau}} P(\mathbf{a}^*_k, \mathbf{a}^*_j)\mathbf{p}_j$, for all $k \in [\boldsymbol{\tau}]$. Then $\mathbb{P}\{\lim_{t\to\infty} \overline{\mathbf{a}}_t = \mathbf{a}^*_k\} = \mathbf{p}_k$, for all $k \in [\boldsymbol{\tau}]$.* □

## APPENDIX C
## PROOF OF PROPOSITION 1: EQUILIBRIUM ACTION PROFILES

Any equilibrium $\overline{a}^* := (a^*_1, \ldots, a^*_n)$ of (7) should satisfy $a^*_i = \text{sign}(\sum_{j\in\mathcal{N}(i)} w_j a^*_j + \eta_i + \lambda_1(\mathbf{s}_{i,t}))$, with probability one for all $i$ and $t$. Hence, $\overline{a}^* \in \{\pm 1\}^n$ is an equilibrium of (7) if, and only if, $-\lambda_1(s_i) \leq \sum_{j\in\mathcal{N}(i)} w_j a^*_j + \eta_i, \forall s_i \in \mathcal{S}_i$ whenever $a^*_i = 1$, and $-\lambda_1(s_i) \geq \sum_{j\in\mathcal{N}(i)} w_j a^*_j + \eta_i, \forall s_i \in \mathcal{S}_i$ whenever $a^*_i = -1$. By multiplying both sides of the inequalities by $a^*_i$ in each case and reordering the terms we derive the claimed characterization of the equilibria or the absorbing states under the action update dynamics in (7).[1] □

## APPENDIX D
## PROOF OF LEMMA 2: TIME-ONE BAYESIAN BELIEFS

We begin by applying the Bayes rule to the observation of agent $i$ at time 1 which include her neighbors' initial beliefs $\{\boldsymbol{\mu}_{j,0}(\cdot); j \in \mathcal{N}(i)\}$ as well as her private signal $\mathbf{s}_{i,0}$. Accordingly, for any $\hat{\theta} \in \Theta$:

$$\boldsymbol{\mu}_{i,1}(\hat{\theta}) = \mathbb{P}_{i,0}\left(\hat{\theta} \mid \mathbf{s}_{i,0}, \{\boldsymbol{\mu}_{j,0}(\cdot); j \in \mathcal{N}(i)\}\right) \quad (22)$$
$$= \frac{\mathbb{P}_{i,0}(\hat{\theta}, \mathbf{s}_{i,0}, \{\boldsymbol{\mu}_{j,0}(\cdot); j \in \mathcal{N}(i)\})}{\mathbb{P}_{i,0}(\mathbf{s}_{i,0}, \{\boldsymbol{\mu}_{j,0}(\cdot); j \in \mathcal{N}(i)\})},$$
$$= \frac{\mathbb{P}_{i,0}(\hat{\theta}, \mathbf{s}_{i,0}, \{\boldsymbol{\mu}_{j,0}(\cdot); j \in \mathcal{N}(i)\})}{\sum_{\tilde{\theta}\in\Theta} \mathbb{P}_{i,0}(\tilde{\theta}, \mathbf{s}_{i,0}, \{\boldsymbol{\mu}_{j,0}(\cdot); j \in \mathcal{N}(i)\})}.$$

---

[1] Here, and in writing the conditions for the case of $a^*_i = -1$ as non-strict inequalities we have violated our earlier convention that agents choose $+1$ when they are indifferent between $+1$ and $-1$. Instead, we are assuming that ties are broken in favor of the equilibrium action profile. This assumption facilitates compact expression of the characterizing conditions for the equilibrium action profiles, and it will have no effect unless with some pathological settings of the signal structure and priors leading to $\sum_{j\in\mathcal{N}(i)} w_j a^*_j + \eta_i = 0$ for some $s_i \in \mathcal{S}_i$, $i \in [n]$.

The succeeding steps follow those in [99] for the case of two communicating agents. For any $j \in [n]$ and all $\pi(\cdot) \in \Delta\Theta$, define the correspondence $\mathcal{I}_j : \Delta\Theta \to \mathscr{P}(\mathcal{S}_j)$ and function $\mathcal{K}_j : \Delta\Theta \to \mathbb{R}$, given by:

$$\mathcal{I}_j(\pi(\cdot)) = \{s \in \mathcal{S}_j : \pi(\hat{\theta}) = \frac{\nu_j(\hat{\theta})l_j(s \mid \hat{\theta})}{\sum_{\tilde{\theta}\in\Theta}\nu_j(\tilde{\theta})l_j(s \mid \tilde{\theta})}, \forall \hat{\theta} \in \Theta\},$$
$$\mathcal{K}_j(\pi(\cdot)) = \sum_{s\in\mathcal{I}_j(\pi(\cdot))} \sum_{\tilde{\theta}\in\Theta} \nu_j(\tilde{\theta})l_j(s \mid \tilde{\theta}). \quad (23)$$

In (23), $\mathcal{I}_j(\pi(\cdot))$ signifies the set of private signals for agent $j$, which are consistent with the observation of belief $\pi(\cdot)$ in that agent. By the same token, $\mathcal{K}_j(\pi(\cdot))$ in (23) is the ex-ante probability for the event that the private signal of agent $j$ belongs to the set $\mathcal{I}_j(\pi(\cdot))$.

The terms $\mathbb{P}_{i,0}(\tilde{\theta}, \mathbf{s}_{i,0}, \{\boldsymbol{\mu}_{j,0}(\cdot); j \in \mathcal{N}(i)\})$ for $\tilde{\theta} \in \Theta$, which appear in the both numerator and denominator of (22) can be simplified by conditioning on the neighbors' observed signals $\{\mathbf{s}_{j,0}; j \in \mathcal{N}(i)\}$ as follows in (24).

$$\mathbb{P}_{i,0}(\theta, \mathbf{s}_{i,0}, \{\boldsymbol{\mu}_{j,0}(\cdot); j \in \mathcal{N}(i)\}) = \sum_{\substack{s_j \in \mathcal{S}_j, \\ j\in\mathcal{N}(i)}}$$
$$\mathbb{P}(\theta, \mathbf{s}_{i,0}, \{\boldsymbol{\mu}_{j,0}(\cdot); j \in \mathcal{N}(i)\} \mid \{\mathbf{s}_{j,0} = s_j; j \in \mathcal{N}(i)\}) \times \ldots$$
$$\ldots \mathbb{P}_{i,0}(\{\mathbf{s}_{j,0} = s_j; j \in \mathcal{N}(i)\}). \quad (24)$$

We next express $\mathbb{P}_{i,0}(\cdot)$ in terms of the priors and signal structures leading to:

$$\mathbb{P}_{i,0}(\tilde{\theta}, \mathbf{s}_{i,0}, \{\boldsymbol{\mu}_{j,0}(\cdot); j \in \mathcal{N}(i)\}) \quad (25)$$
$$= \sum_{\{s_j \in \mathcal{I}_j(\boldsymbol{\mu}_{j,0}(\cdot)), j\in\mathcal{N}(i)\}} \nu_i(\tilde{\theta})l_i(\mathbf{s}_{i,0} \mid \tilde{\theta}) \prod_{j\in\mathcal{N}(i)} l_j(s_j \mid \tilde{\theta})$$
$$= \frac{\nu_i(\tilde{\theta})l_i(\mathbf{s}_{i,0} \mid \tilde{\theta})}{\prod_{j\in\mathcal{N}(i)} \nu_j(\tilde{\theta})} \prod_{j\in\mathcal{N}(i)} \left( \sum_{\substack{s_j \in \\ \mathcal{I}_j(\boldsymbol{\mu}_{j,0}(\cdot))}} \nu_j(\tilde{\theta})l_j(s_j \mid \tilde{\theta}) \right).$$

Bayes rule in (1), together with the functions defined in (23), can now be used to eliminate the product terms involving $s_j$ from (25) and get:

$$\mathbb{P}_{i,0}(\tilde{\theta}, \mathbf{s}_{i,0}, \{\boldsymbol{\mu}_{j,0}(\cdot); j \in \mathcal{N}(i)\}) = \frac{\nu_i(\tilde{\theta})l_i(\mathbf{s}_{i,0} \mid \tilde{\theta})}{\prod_{j \in \mathcal{N}(i)} \nu_j(\tilde{\theta})} \times \quad (26)$$

$$\ldots \times \prod_{j \in \mathcal{N}(i)} \left( \boldsymbol{\mu}_{j,0}(\tilde{\theta}) \sum_{s_j \in \mathcal{I}_j(\boldsymbol{\mu}_{j,0}(\cdot))} \sum_{\bar{\theta} \in \Theta} \nu_j(\bar{\theta})l_j(s \mid \bar{\theta}) \right)$$

$$= \nu_i(\tilde{\theta})l_i(\mathbf{s}_{i,0} \mid \tilde{\theta}) \left( \prod_{j \in \mathcal{N}(i)} \frac{\boldsymbol{\mu}_{j,0}(\tilde{\theta})}{\nu_j(\tilde{\theta})} \right) \prod_{j \in \mathcal{N}(i)} \mathcal{K}_j(\boldsymbol{\mu}_{j,0}(\cdot)).$$

Upon replacing (26) in (22), the product terms involving $\mathcal{K}_j(\boldsymbol{\mu}_{j,0}(\cdot))$ cancel out and (9) follows. □

## APPENDIX E
## PROOF OF THEOREM 1: NO LEARNING WHEN $\rho > 1$

We begin the analysis of the beliefs propagation under (10) by forming the ratio

$$\frac{\boldsymbol{\mu}_{i,t}(\check{\theta})}{\boldsymbol{\mu}_{i,t}(\theta)} = \frac{\nu_i(\check{\theta})}{\nu_i(\theta)} \times \frac{l_i(\mathbf{s}_{i,t} \mid \check{\theta})}{l_i(\mathbf{s}_{i,t} \mid \theta)} \times \prod_{j \in \mathcal{N}(i)} \frac{\boldsymbol{\mu}_{j,t-1}(\check{\theta})}{\boldsymbol{\mu}_{j,t-1}(\theta)} \times \frac{\nu_j(\theta)}{\nu_j(\check{\theta})},$$

for any false state $\check{\theta} \in \Theta \setminus \{\theta\}$ and each agent $i \in [n]$ at all times $t \in \mathbb{N}$. The above has the advantage of removing the normalization factor in the dominator out of the picture; thence, focusing instead on the evolution of belief ratios, which has a log-linear format. The latter motivates definitions of log-likelihood ratios for signals, beliefs, and priors as follows. Similarly to $\lambda_{\check{\theta}}(\mathbf{s}_{i,t})$ and $\phi_{i,t}(\check{\theta})$, define the log-ratios of prior beliefs as $\gamma_i(\check{\theta}) := \log(\nu_i(\check{\theta})/\nu_i(\theta))$. Starting from the above iterations for the belief ratio and taking the logarithms of both sides yields

$$\phi_{i,t}(\check{\theta}) = \gamma_i(\check{\theta}) + \lambda_{\check{\theta}}(\mathbf{s}_{i,t}) + \sum_{j \in \mathcal{N}(i)} \phi_{j,t-1}(\check{\theta}) - \gamma_j(\check{\theta}). \quad (27)$$

Multiplying both sides of (27) by $\alpha_i$, which is the centrality of agent $i$, and summing over all $i \in [n]$ yields that

$$\boldsymbol{\Phi}_t(\check{\theta}) = \sum_{i=1}^n \alpha_i \gamma_i(\check{\theta}) + \sum_{i=1}^n \alpha_i \lambda_{\check{\theta}}(\mathbf{s}_{i,t}) \quad (28)$$

$$+ \sum_{i=1}^n \alpha_i \sum_{j \in \mathcal{N}(i)} (\phi_{j,t-1}(\check{\theta}) - \gamma_j(\check{\theta})).$$

First note that we can write

$$\sum_{i=1}^n \alpha_i \gamma_i(\check{\theta}) - \sum_{i=1}^n \alpha_i \sum_{j \in \mathcal{N}(i)} \gamma_j(\check{\theta}) = \text{tr}\left\{ \left(I - A^T\right) \overline{\alpha} \overline{\gamma}(\check{\theta})^T \right\}$$

$$= (1-\rho)\beta(\check{\theta}), \quad (29)$$

where $\overline{\gamma}(\check{\theta}) := (\gamma_1(\check{\theta}), \ldots, \gamma_n(\check{\theta}))^T$. Next note that by the choice of $\overline{\alpha}$ as the eigenvector corresponding to the $\rho$ eigen-

value of matrix $A$ we get

$$\sum_{i=1}^n \alpha_i \sum_{j \in \mathcal{N}(i)} \phi_{j,t-1}(\check{\theta}) = \overline{\alpha}^T A \overline{\phi}_{t-1}(\check{\theta}) = \rho \overline{\alpha}^T \overline{\phi}_{t-1}(\check{\theta})$$

$$= \rho \boldsymbol{\Phi}_{t-1}(\check{\theta}). \quad (30)$$

where $\overline{\phi}_t(\check{\theta}) := (\phi_{1,t}(\check{\theta}), \ldots, \phi_{n,t}(\check{\theta}))^T$. Now replacing (29) and (30) in (28) yields the following recursion for $\boldsymbol{\Phi}_t(\check{\theta})$:

$$\boldsymbol{\Phi}_t(\check{\theta}) = \boldsymbol{\Lambda}_t(\check{\theta}) + \rho \boldsymbol{\Phi}_{t-1}(\check{\theta}) + (1-\rho)\beta(\check{\theta}), \quad (31)$$

initialized by $\boldsymbol{\Phi}_0(\check{\theta}) = \beta(\check{\theta}) + \boldsymbol{\Lambda}_0(\check{\theta})$, where $\beta(\check{\theta}) := \sum_{i=1}^n \alpha_i \log(\nu_i(\check{\theta})/\nu_i(\theta))$ is a constant that is determined by the initial prior beliefs, and it measures the total bias in the network relative between the two states $\check{\theta}$ and $\theta$. In particular, if the agents are unbiased starting from uniform priors on $\Theta$, then $\beta(\check{\theta}) = 0, \forall \check{\theta} \in \Theta$. Note also that the assumption of full support priors implies that $|\beta(\check{\theta})|$ is finite. By iterating (31) for $t \in \mathbb{N}$ we obtain (14). Next note that in a strongly connected graph every node has a degree greater than or equal to one so that $\rho \geq 1$, [84, Chapter 2]. If $\rho > 1$, then the term $\rho^t \boldsymbol{\Lambda}_0(\check{\theta})$ increases in variance as $t \to \infty$, and unless $\boldsymbol{\Lambda}_0(\check{\theta}) < \epsilon$ with $\mathbb{P}_\theta$-probability one for some $\epsilon < 0$, almost sure convergence to $-\infty$ for $\boldsymbol{\Phi}_t(\check{\theta})$ in (14) cannot hold true. □